\documentclass[reqno,12pt]{amsart}

\newcommand{\hyp}[5]{\,\mbox{}_{#1}F_{#2}\!\left(
  \genfrac{}{}{0pt}{}{#3}{#4};#5\right)}

\newcommand{\Ga}{\Gamma}

\newcommand{\ta}{\theta}
\newcommand{\pa}{\partial}
\newcommand{\rd}{\,\mathrm d}

\usepackage{amsmath,amssymb}

\allowdisplaybreaks

\title[Multiple hypergeometric series]{Multiple hypergeometric series
-- Appell series and beyond}
\author{Michael J.\ Schlosser}
\address{Fakult\"at f\"ur Mathematik,
Universit\"at Wien, Nordbergstrasse 15, A-1090 Vienna, Austria}
\email{michael.schlosser@univie.ac.at}
\thanks{Partially supported by FWF Austrian Science Fund
grants S9607 \& F50-08.}

\begin{document}

\begin{abstract}
This survey article (which will appear as a chapter in the book
``Computer Algebra in Quantum Field Theory: Integration, Summation and
Special Functions'', Springer-Verlag) provides a small collection
of basic material on multiple hypergeometric series of Appell-type
and of more general series of related type.
\end{abstract}

\maketitle

\section{Introduction}\label{sec:1}

Hypergeometric series and its various generalizations,
in particular such involving \textit{multiple} series, appear
in various branches of mathematics and its applications.
This survey article features a small collection of selected material
on multiple hypergeometric series of \textit{Appell}-type and
of more general series of closely related type.

These types of series appear very naturally in quantum
field theory, in particular in
the computation of analytic expressions for Feynman integrals
(for which we kindly refer to other relevent chapters in this volume).
Such integrals can be obtained and computed in different ways
-- which may lead to identities for Appell series
(see e.g.~M.A.~Shpot~\cite{SchlosserShpot:07}).
 On the other hand, the
application of known relations for Appell series may lead
to simplifications, help to solve problems or lead to more
insight in quantum field theory. Therefore it is of importance
that people working in this area have a basic understanding of the
existing theory for such series\footnote{Researchers
working with Feynman integrals who are in demand of effective
manipulation of Appell-type series including
differential reductions and $\epsilon$-expansions
may find HYPERDIRE (located at
\texttt{https://sites.google.com/site/loopcalculations/})
useful, which is a set of Wolfram Mathematica
based programs
for differential reduction of Horn-type hypergeometric functions,
see V.~Bytev et al.~\cite{SchlosserBytev+al:11}.}.
This survey is meant to 
provide a very digestible, easy introduction to Appell-type series.
Besides of recalling some known results including various useful
identities satisfied by the series, some of the standard
mathematical techniques which are used to prove and derive
these identities are illustrated.
We highlight some of the most fundamental
properties and relations for Appell hypergeometric series and
further give hints of similar relations for the series which are
(slightly) beyond the hierarchy of Appell series.
All the series we consider admit very explicit series and integral
representations.

To warn the reader: There exist various different types of multivariate
hypergeometric series which are not covered in this survey.
In particular, here we do not treat multiple hypergeometric series
associated with \textit{root systems}~\cite{SchlosserHeckmanOpdam:87,SchlosserKoornwinder:98,SchlosserMilne:01,SchlosserRosengren:04}, hypergeometric series of
{\em matrix argument}~\cite{SchlosserGrossRichards:89},
and other types of multivariate hypergeometric series such
as those which mainly appear in the study of {\em orthogonal polynomials}
of severable variables (often also associated with root systems)
\cite{SchlosserDunklXu:01,SchlosserMacdonald:98}.

A very important extension of Appell-type series which is just
beyond the scope of this basic survey article are the
multivariate hypergeometric functions considered by
I.M.~Gelfand, M.M.~Kapranov, and
A.V.~Zelevinsky~\cite{SchlosserGKZ:89}, developed in the
late 1980's. These $A$-hypergeometric functions
(or GKZ-hypergeometric functions) are fundamental objects
in the theory of integrable systems as they are the holonomic
solutions of a (certain) $A$-hypergeometric system of
partial differential equations.
Natural questions regarding algebraic solutions and
monodromy for $A$-hypergeometric functions have been
recently addressed by
F.~Beukers~\cite{SchlosserBeukers:10,SchlosserBeukers:11}.

For {\em basic} (or $q$-{\em series}) analogues of Appell functions,
see G.~Gasper and M.~Rahman's text~\cite[Ch.~10]{SchlosserGasperRahman:04}.

\section{Appell series}\label{sec:2}
Appell series are a natural two-variable extension of
hypergeometric series. They are treated with detail in \'Erdelyi et
al.~\cite{SchlosserErdelyi+al:81}, the classical reference
for special functions. 

In the following, we follow to great extent the expositions from
the classical texts of W.N.~Bailey~\cite{SchlosserBailey:64}, and
L.J.~Slater~\cite{SchlosserSlater:66} (both contain a great amount
of material on hypergeometric series).

For convenience, we use the Pochhammer symbol notation for the
shifted factorial,

\begin{subequations}
\begin{equation}
(a)_n:=\begin{cases}
a(a+1)\ldots(a+n-1)\quad &\mbox{if $n=1,2,\ldots\,$,}\\
1&\mbox{if $n=0$.}
\end{cases}
\end{equation}
Accordingly, we have
\begin{equation}
(a)_n=\frac{\Gamma(a+n)}{\Gamma(a)}
\end{equation}
\end{subequations}
which is used as a definition for the shifted factorial in case $n$ is not
necessarily a nonnegative integer.


The goal is to generalize the Gau{\ss} hypergeometric function
$$
\hyp 21{a,\ b}cx=\sum_{n\ge0}\frac{(a)_n\,(b)_n}{n!\,(c)_n}x^n
$$
to a double series depending on two variables.

The easiest is to consider the simple product
$$
\hyp 21{a,\ b}cx \hyp 21{a',\ b'}{c'}y=
\sum_{m\ge0}\sum_{n\ge0}
\frac{(a)_m\,(a')_n\,(b)_m\,(b')_n}{m!\,n!\ (c)_m\,(c')_n}x^my^n,
$$
where on the right-hand side the indices $m,n$ appear
uncoupled.

To consider a genuine double series instead
(which does not factor into a simple product of two series),
we now deliberately choose to replace one or more of the three
products $(a)_m\,(a')_n$, $(b)_m\,(b')_n$, $(c)_m\,(c')_n$
by products of coupled type $(a)_{m+n}$ (other choices
such as $(a)_{m-n}$ or $(a)_{2m-n}$, etc., instead,
may be sensible as well;
they lead to Horn-type series, see Subsection~\ref{subsec:Horn}). 

There are five different possibilities,
one of which by application of the binomial
theorem gives the series
$$
\sum_{m\ge0}\sum_{n\ge0}
\frac{(a)_{m+n}\,(b)_{m+m}}{m!\,n!\ (c)_{m+n}}x^my^n=
\hyp 21{a,\ b}c{x+y},
$$
i.e., an ordinary hypergeometric series.

The other four remaining possibilities are classified as $F_1$-,
$F_2$-, $F_3$-, and $F_4$-series
(cf.\ P.~Appell~\cite{SchlosserAppell:80} and P.~Appell \&
M.-J.~Kamp\'e de F\'eriet~\cite{SchlosserAppellKdF:26}):

\begin{subequations}
\begin{alignat}2
F_1\big(a;b,b';c;x,y\big)&:=\sum_{m\ge0}\sum_{n\ge0}
\frac{(a)_{m+n}\,(b)_m\,(b')_n}{m!\,n!\ (c)_{m+n}}x^my^n,
& |x|,|y|<1.\\
F_2\big(a;b,b';c,c';x,y\big)&:=\sum_{m\ge0}\sum_{n\ge0}
\frac{(a)_{m+n}\,(b)_m\,(b')_n}{m!\,n!\ (c)_m\,(c')_n}x^my^n,
& |x|+|y|<1.\\
F_3\big(a,a';b,b';c;x,y\big)&:=\sum_{m\ge0}\sum_{n\ge0}
\frac{(a)_m\,(a')_n\,(b)_m\,(b')_n}{m!\,n!\ (c)_{m+n}}x^my^n,
& |x|,|y|<1.\\
F_4\big(a;b;c,c';x,y\big)&:=\sum_{m\ge0}\sum_{n\ge0}
\frac{(a)_{m+n}\,(b)_{m+n}}{m!\,n!\ (c)_m\,(c')_n}x^my^n,
& |x|^{\frac 12}+|y|^{\frac 12}<1.
\end{alignat}
\end{subequations}

One immediately observes the following simple identities:
\begin{equation}
F_1\big(a;b,b';c;x,y\big)=\sum_{m\ge 0}\frac{(a)_m\,(b)_m}{m!\,(c)_m}x^m
\hyp 21{a+m,\,b'}{c+m}y.
\end{equation}
\begin{subequations}
\begin{align}
F_1\big(a;b,b';c;x,0\big)=F_2\big(a;b,b';c,c';x,0\big)=
F_3\big(a,a';b,b';c;x,0\big)&\\
=F_4\big(a;b;c,c';x,0\big)=\hyp 21{a,\,b}cx&.
\end{align}
\end{subequations}
\begin{subequations}
\begin{align}
F_1\big(a;b,0;c;x,y\big)=F_2\big(a;b,0;c,c';x,y\big)
=F_3\big(a,a';b,0;c;x,y\big)&\\
=\hyp 21{a,\,b}cx&.
\end{align}
\end{subequations}

Using ideas of N.Ja.~Vilenkin~\cite{SchlosserVilenkin:68},
W.~Miller, Jr.~\cite{SchlosserMiller:73} has given a Lie theoretic
interpretation of the Appell functions $F_1$. In particular,
he showed that $sl(5,\mathbb C)$ is the dynamical symmetry algebra
for the $F_1$.

\subsection{Contiguous relations and recursions}
All contiguous relations
for the $F_1$ function can be derived from these four relations:
\begin{subequations}
\begin{align}
(a-b-b')\, F_1\big(a;b,b';c;x,y\big)-a\, F_1\big(a+1;b,b';c;x,y\big)&
\notag\\
+b\, F_1\big(a;b+1,b';c;x,y\big)+b'\, F_1\big(a;b,b'+1;c;x,y\big)&=0,\\
c\, F_1\big(a;b,b';c;x,y\big)-(c-a)\, F_1\big(a;b,b';c+1;x,y\big)&
\notag\\
-a\, F_1\big(a+1;b,b';c+1;x,y\big)&=0,\\
c\, F_1\big(a;b,b';c;x,y\big)+c(x-1)\, F_1\big(a;b+1,b';c;x,y\big)&
\notag\\
-(c-a)x\, F_1\big(a;b+1,b';c+1;x,y\big)&=0,\\
c\, F_1\big(a;b,b';c;x,y\big)+c(y-1)\, F_1\big(a;b,b'+1;c;x,y\big)&
\notag\\
-(c-a)y\, F_1\big(a;b,b'+1;c+1;x,y\big)&=0.
\end{align}
\end{subequations}
Similar sets of relations exist for the other Appell functions, see
R.G.~Buschman~\cite{SchlosserBuschman:87}.

Recently, X.~Wang~\cite{SchlosserWang:12} has used contiguous relations
and induction to derive various recursion formulae
for all the Appell
functions $F_1,F_2,F_3,F_4$. (Some of the recursions for $F_2$ were
previously given by S.B.~Opps, N.~Saad and
H.M.~Srivastava~\cite{SchlosserOpps+al:09}.)
For $n=1$ these recursions reduce to equivalent forms of the
known contiguous relations.

In particular, for $F_1$ we have
\begin{subequations}
\begin{align}
F_1\big(a+n;b,b';c;x,y\big)={}&F_1\big(a;b,b';c;x,y\big)+
\frac{bx}c\sum_{k=1}^nF_1\big(a+k;b+1,b';c+1;x,y\big)\notag\\&+
\frac{b'y}c\sum_{k=1}^nF_1\big(a+k;b,b'+1;c+1;x,y\big),\\
F_1\big(a-n;b,b';c;x,y\big)={}&F_1\big(a;b,b';c;x,y\big)-
\frac{bx}c\sum_{k=1}^{n-1}F_1\big(a-k;b+1,b';c+1;x,y\big)\notag\\&-
\frac{b'y}c\sum_{k=1}^{n-1}F_1\big(a-k;b,b'+1;c+1;x,y\big),\\
F_1\big(a;b+n,b';c;x,y\big)={}&F_1\big(a;b,b';c;x,y\big)+
\frac{ax}c\sum_{k=1}^nF_1\big(a+1;b+k,b';c+1;x,y\big),\\
F_1\big(a;b-n,b';c;x,y\big)={}&F_1\big(a;b,b';c;x,y\big)-
\frac{ax}c\sum_{k=1}^{n-1}F_1\big(a+1;b-k,b';c+1;x,y\big),\\
F_1\big(a;b,b';c-n;x,y\big)={}&F_1\big(a;b,b';c;x,y\big)\notag\\&+
abx\sum_{k=1}^n\frac{F_1\big(a+1;b+1,b';c-k+2;x,y\big)}{(c-k)(c-k+1)}
\notag\\&+
ab'y\sum_{k=1}^n\frac{F_1\big(a+1;b,b'+1;c-k+2;x,y\big)}{(c-k)(c-k+1)}.
\end{align}
\end{subequations}

For $F_2$ we have
\begin{subequations}
\begin{align}
F_2\big(a+n;b,b';c,c';x,y\big)={}&F_2\big(a;b,b';c,c';x,y\big)\notag\\&+
\frac{bx}c\sum_{k=1}^nF_2\big(a+k;b+1,b';c+1,c';x,y\big)\notag\\&+
\frac{b'y}{c'}\sum_{k=1}^nF_2\big(a+k;b,b'+1;c,c'+1;x,y\big),\\
F_2\big(a-n;b,b';c,c';x,y\big)={}&F_2\big(a;b,b';c,c';x,y\big)\notag\\&-
\frac{bx}c\sum_{k=1}^{n-1}F_2\big(a-k;b+1,b';c+1,c';x,y\big)\notag\\&-
\frac{b'y}{c'}\sum_{k=1}^{n-1}F_2\big(a+k;b,b'+1;c,c'+1;x,y\big),\\
F_2\big(a;b+n,b';c,c';x,y\big)={}&F_2\big(a;b,b';c,c';x,y\big)\notag\\&+
\frac{ax}c\sum_{k=1}^nF_2\big(a+1;b+k,b';c+1,c';x,y\big),\\
F_2\big(a;b-n,b';c,c';x,y\big)={}&F_2\big(a;b,b';c,c';x,y\big)\notag\\&-
\frac{ax}c\sum_{k=1}^{n-1}F_2\big(a+1;b-k,b';c+1,c';x,y\big),\\
F_2\big(a;b,b';c-n,c';x,y\big)={}&F_2\big(a;b,b';c,c';x,y\big)\notag\\&+
abx\sum_{k=1}^n\frac{F_2\big(a+1;b+1,b';c-k+2,c';x,y\big)}{(c-k)(c-k+1)}.
\end{align}
\end{subequations}

For $F_3$ we have
\begin{subequations}
\begin{align}
F_3\big(a+n,a';b,b';c;x,y\big)={}&F_3\big(a,a';b,b';c;x,y\big)\notag\\&+
\frac{bx}c\sum_{k=1}^nF_3\big(a+k,a';b+1,b';c+1;x,y\big),\\
F_3\big(a-n,a';b,b';c;x,y\big)={}&F_3\big(a,a';b,b';c;x,y\big)\notag\\&-
\frac{bx}c\sum_{k=1}^{n-1}F_3\big(a-k,a';b+1,b';c+1;x,y\big),\\
F_3\big(a,a';b,b';c-n;x,y\big)={}&F_3\big(a,a';b,b';c;x,y\big)\notag\\&+
abx\sum_{k=1}^n\frac{F_3\big(a+1,a';b+1,b';c-k+2;x,y\big)}{(c-k)(c-k+1)}
\notag\\&+
a'b'y\sum_{k=1}^n\frac{F_3\big(a,a'+1;b,b'+1;c-k+2;x,y\big)}{(c-k)(c-k+1)}.
\end{align}
\end{subequations}

Finally, for $F_4$ we have
\begin{subequations}
\begin{align}
F_4\big(a+n;b;c,c';x,y\big)={}&F_4\big(a;b;c,c';x,y\big)\notag\\&+
\frac{bx}c\sum_{k=1}^nF_4\big(a+k;b+1;c+1,c';x,y\big)\notag\\&+
\frac{by}{c'}\sum_{k=1}^nF_4\big(a+k;b+1;c,c'+1;x,y\big),\\
F_4\big(a-n;b;c,c';x,y\big)={}&F_4\big(a;b;c,c';x,y\big)\notag\\&-
\frac{bx}c\sum_{k=1}^{n-1}F_4\big(a-k;b+1;c+1,c';x,y\big)\notag\\&-
\frac{by}{c'}\sum_{k=1}^{n-1}F_4\big(a-k;b+1;c,c'+1;x,y\big),\\
F_4\big(a;b;c-n,c';x,y\big)={}&F_4\big(a;b;c,c';x,y\big)\notag\\&+
abx\sum_{k=1}^{n-1}\frac{F_4\big(a+1;b+1;c-k+1,c';x,y\big)}{(c-k)(c-k-1)}.
\end{align}
\end{subequations}

Most of these recursions can be extended to elegant recursions
involving more terms.
For instance,
\begin{subequations}
\begin{align}
F_1\big(a+n;b,b';c;x,y\big)={}&\sum_{i=0}^n\sum_{k=0}^{n-i}\binom ni
\binom{n-i}k\frac{(b)_i(b')_k}{(c)_{k+i}}\notag\\&\times x^iy^j
F_1\big(a+i+k;b+i,b'+k;c+i+k;x,y\big),\\
F_1\big(a-n;b,b';c;x,y\big)={}&\sum_{i=0}^n\sum_{k=0}^{n-i}\binom ni
\binom{n-i}k\frac{(b)_i(b')_k}{(c)_{k+i}}\notag\\&\times (-x)^i(-y)^j
F_1\big(a;b+i,b'+k;c+i+k;x,y\big),\\
F_1\big(a;b+n,b';c;x,y\big)={}&\sum_{k=0}^n\binom nk
\frac{(a)_k}{(c)_k}x^kF_1\big(a+k;b+k,b';c+k;x,y\big),\\
F_1\big(a;b-n,b';c;x,y\big)={}&\sum_{k=0}^n\binom nk
\frac{(a)_k}{(c)_k}(-x)^kF_1\big(a+k;b,b';c+k;x,y\big),
\end{align}
\end{subequations}
or
\begin{equation}
F_4\big(a;b;c-n,c';x,y\big)=\sum_{k=0}^n\binom nk
\frac{(a)_k(b)_k}{(c)_k(c-n)_k}x^k
F_4\big(a+k;b+k;c+k,c';x,y\big).
\end{equation}

\subsection{Partial differential equations}
Let
$$
z=F_1\big(a;b,b';c;x,y\big)=\sum_{m\ge0}\sum_{n\ge0}A_{m,n}x^my^n.
$$
Then
$$
A_{m+1,n}=\frac{(a+m+n)(b+m)}{(1+m)(c+m+n)}A_{m,n},
$$
and
$$
A_{m,n+1}=\frac{(a+m+n)(b'+n)}{(1+n)(c+m+n)}A_{m,n}.
$$

Denoting the partial differential operators by
$$\ta=x\frac{\pa}{\pa x}\qquad\text{and}
\qquad\phi=y\frac{\pa}{\pa y},$$
we readily see that $z=F_1$ satisfies the
partial differential equations
\begin{subequations}
\begin{align}
\big[(\ta+\phi+a)(\ta+b)-\frac 1x\ta(\ta+\phi+c-1)\big]z&=0,\\
\big[(\ta+\phi+a)(\phi+b')-\frac 1x\phi(\ta+\phi+c-1)\big]z&=0.
\end{align}
\end{subequations}

Now let
$$p=\frac{\pa z}{\pa x},\quad q=\frac{\pa z}{\pa y},\quad
r=\frac{\pa z}{\pa x}\frac{\pa z}{\pa y},\quad
s=\frac{\pa z^2}{\pa x^2},\quad t=\frac{\pa z^2}{\pa y^2}.
$$

Then $z=F_1$ satisfies the partial differential equations
\begin{subequations}
\begin{align}
x(1-x)r+y(1-x)s+[c-(a+b+1)x]p-byq-abz&=0,\\
y(1-y)t+x(1-y)s+[c-(a+b'+1)y]q-b'xp-ab'z&=0.
\end{align}
\end{subequations}

Similarly, $z=F_2$ satisfies the partial differential equations
\begin{subequations}
\begin{align}
x(1-x)r-xys+[c-(a+b+1)x]p-byq-abz&=0,\\
y(1-y)t-xys+[c'-(a+b'+1)y]q-b'xp-ab'z&=0.
\end{align}
\end{subequations}

Similarly, $z=F_3$ satisfies the partial differential equations
\begin{subequations}
\begin{align}
x(1-x)r+ys+[c-(a+b+1)x]p-abz&=0,\\
y(1-y)t+xs+[c-(a'+b'+1)y]q-a'b'z&=0.
\end{align} 
\end{subequations}

Finally, $z=F_4$ satisfies the partial differential equations
\begin{subequations}
\begin{align}
x(1-x)r-y^2t-2xys+cp-(a+b+1)(xp+yq)-abz&=0,\\
y(1-y)t-x^2r-2xys+c'q-(a+b+1)(xp+yq)-abz&=0.
\end{align} 
\end{subequations}

\subsection{Integral representations}
Integral representations for Appell series are very useful.
Substitution of variables in theses integrals lead to equivalent
integrals. This provides an effective and easy method to derive
transformation formulae for Appell series, see Subsection~\ref{subsec:tf}.

Consider the integral
$$
I=\int\!\!\!\int\! u^{b-1}v^{b'-1}(1-u-v)^{c-b-b'-1}(1-ux-vy)^{-a}\rd u\rd v,
$$
taken over the triangular region $u\ge 0$,\, $v\ge 0$,\, $u+v\le1$.
(We implicitly assume suitable conditions of the parameters $a,b,b',c$
such that the integral is well-defined and converges.)

Now, provided $|vy/(1-ux)|<1$, we have, by binomial expansion,
\begin{align*}
(1-ux-vy)^{-a}&=(1-ux)^{-a}\sum_{m\ge0}\frac{(a)_m}{m!}
\left(\frac{vy}{1-ux}\right)^m\\
&=\sum_{m\ge0}\frac{(a)_m}{m!}v^my^m(1-ux)^{-a-m}\\
&=\sum_{m\ge0}\frac{(a)_m}{m!}v^my^m\sum_{n\ge0}\frac{(a+m)_n}{n!}u^nx^n.
\end{align*}

Thus,
\begin{align*}
I&=\sum_{m\ge0}\sum_{n\ge0}\frac{(a)_{m+n}}{m!n!}x^ny^m
\int\!\!\!\int\! u^{b-1+n}v^{b'-1+m}(1-u-v)^{c-b-b'-1}\rd u\rd v\\
&=\sum_{m\ge0}\sum_{n\ge0}\frac{(a)_{m+n}}{m!n!}x^ny^m\,
\Ga\!\begin{bmatrix}b+n,b'+m,c-b-b'\\c+m+n\end{bmatrix},
\end{align*}
which yields
\begin{equation}
I=\Ga\!\begin{bmatrix}b,b',c-b-b'\\c\end{bmatrix}\,F_1\big(a;b,b';c;x,y\big).
\end{equation}
While $I$ is a double integral, a {\em single integral} for $F_1$
even exists, see \eqref{singleintegralF1}.

Similarly,
\begin{align}
\int_0^1\!\!\!\int_0^1\! u^{b-1}v^{b'-1}(1-u)^{c-b'-1}(1-v)^{c'-b'-1}
(1-ux-vy)^{-a}\rd u\rd v&\notag\\
=\Ga\!\begin{bmatrix}b,b',c-b,c'-b'\\c,\,c'\end{bmatrix}\,
F_2\big(a;b,b';c,c';x,y\big)&,
\end{align}
and
\begin{align}
\int\!\!\!\int\! u^{b-1}v^{b'-1}(1-u-v)^{c-b-b'-1}
(1-ux)^{-a}(1-vy)^{-a'}\rd u\rd v&\notag\\
=\Ga\!\begin{bmatrix}b,b',c-b-b'\\c'\end{bmatrix}\,
F_3\big(a,a';b,b';c';x,y\big)&,
\end{align}
the  last integral taken over the triangular region
$u\ge 0$, $v\ge 0$, $u+v\le1$.

The double integral for $F_4$ is more complicated:
\begin{align}
\int_0^1\!\!\!\int_0^1\! u^{a-1}v^{b-1}(1-u)^{c-a-1}
(1-v)^{c'-b-1}(1-ux)^{-b}(1-vy)^{-a}&\notag\\\times
\left(1-\frac{uvxy}{(1-ux)(1-vy)}\right)^{c+c'-a-b-1}
\rd u\rd v&\notag\\
=\Ga\!\begin{bmatrix}a,b,c-a,c'-b\\c,\,c'\end{bmatrix}\,
F_4\big(a;b;c,c';x(1-y),y(1-x)\big)&.
\end{align}

In 1881, \'E.~Picard~\cite{SchlosserPicard:81}
discovered a single integral for $F_1$.
Let
$$
I'=\int_0^1 u^{a-1}(1-u)^{c-a-1}(1-ux)^{-b}(1-uy)^{-b'}\rd u,
$$
where $\Re c>\Re a>0$.
Then
\begin{align*}
I'&=\sum_{m\ge 0}\sum_{n\ge 0}\int_0^1u^{a-1}(1-u)^{c-a-1}
\frac{(b)_m}{m!}u^mx^m\frac{(b')_n}{n!}
u^ny^n\rd u\\
&=\sum_{m\ge 0}\sum_{n\ge 0}\frac{(b)_m(b')_n}{m!n!}x^my^n
\int_0^1u^{a+m+n-1}(1-u)^{c-a-1}\rd u\\
&=\sum_{m\ge 0}\sum_{n\ge 0}\frac{(b)_m(b')_n}{m!n!}x^my^n\,
\Ga\!\begin{bmatrix}a+m+n,c-a\\c+m+n\end{bmatrix},
\end{align*}
hence
\begin{equation}\label{singleintegralF1}
I'=\Ga\!\begin{bmatrix}a,c-a\\c\end{bmatrix}F_1\big(a;b,b';c;x,y\big).
\end{equation}

\subsubsection{Incomplete elliptic integrals}
As immediate consequences of \eqref{singleintegralF1},
it follows that the incomplete elliptic integrals $F$ and $E$
and the complete elliptic integral $\Pi$ can all be expressed
in terms of special cases of the Appell $F_1$ function:
\begin{subequations}
\begin{align}
F(\phi,k):&=\int_0^\phi\frac{\rd\ta}{\sqrt{1-k^2\sin^2\ta}}\notag\\
&=\sin\phi\;F_1\!\left(\frac 12;\frac 12,\frac 12;\frac 32;
\sin^2\phi,k^2\sin^2\phi\right),\qquad |\Re\phi|<\frac{\pi}2,\\
E(\phi,k):&=\int_0^\phi\sqrt{1-k^2\sin^2\ta}\rd\ta\notag\\
&=\sin\phi\;F_1\!\left(\frac 12;\frac 12,-\frac 12;\frac 32;
\sin^2\phi,k^2\sin^2\phi\right),\qquad |\Re\phi|<\frac{\pi}2,\\
\Pi(n,k):&=\int_0^{\pi/2}\frac{\rd\ta}
{(1-n\sin^2\ta)\sqrt{1-k^2\sin^2\ta}}
=\frac{\pi}2\,F_1\!\left(\frac 12;1,\frac 12;1;
n,k^2\right).
\end{align}
\end{subequations}

\subsection{Transformations}\label{subsec:tf}
In the single integral for the $F_1$ series,
\begin{equation*}
F_1\big(a;b,b';c;x,y\big)=\Ga\!\begin{bmatrix}c\\a,c-a\end{bmatrix}
\int_0^1 u^{a-1}(1-u)^{c-a-1}(1-ux)^{-b}(1-uy)^{-b'}\rd u,
\end{equation*}
one may use the substitution of variables $u=1-v$
to prove
\begin{equation}
F_1\big(a;b,b';c;x,y\big)=(1-x)^{-b}(1-y)^{-b'}
F_1\!\left(c-a;b,b';c;\frac x{x-1},\frac y{y-1}\right)\!.
\end{equation}
For $b'=0$ this reduces to the well-known {\em Pfaff--Kummer transformation}
for the $_2F_1$:
$$
\hyp 21{a,\,b}cx=(1-x)^{-b}\hyp 21{c-a,\,b}c{\frac x{x-1}}.
$$

Similarly, the substitution of variables $u=\frac v{1-x+vx}$
can be used to prove
\begin{equation}
F_1\big(a;b,b';c;x,y\big)=(1-x)^{-a}
F_1\!\left(a;-b-b'+c,b';c;\frac x{x-1},\frac{y-x}{1-x}\right).
\end{equation}
For $b'=0$ this reduces again to the
Pfaff--Kummer transformation for the $_2F_1$ series.

On the other hand, if $c=b+b'$, then
\begin{subequations}
\begin{align}
F_1\big(a;b,b';b+b';x,y\big)&=(1-x)^{-a}
\hyp 21{a,\,b'}{b+b'}{\frac{y-x}{1-x}}\\
&=(1-y)^{-a}\hyp 21{a,\,b}{b+b'}{\frac{x-y}{1-y}}.
\end{align}
\end{subequations}

Similarly,
\begin{equation}
F_1\big(a;b,b';c;x,y\big)=(1-y)^{-a}
F_1\!\left(a;b,c-b-b';c;\frac{x-y}{1-y},\frac y{y-1}\right),
\end{equation}
\begin{equation}
F_1\big(a;b,b';c;x,y\big)={\mbox{\small $(1-x)^{c-a-b}(1-y)^{-b'}$}}
F_1\!\left(c-a;c-b-b',b';c;x,\frac{x-y}{1-y}\right),
\end{equation}
\begin{equation}
F_1\big(a;b,b';c;x,y\big)={\mbox{\small $(1-x)^{-b}(1-y)^{c-a-b'}$}}
F_1\!\left(c-a;b,c-b-b';c;\frac{y-x}{1-x},y\right).
\end{equation}

Further,
\begin{equation}
F_2\big(a;b,b';c,c';x,y\big)=(1-x)^{-a}
F_2\!\left(a;c-b,b';c,c';\frac x{x-1},\frac y{1-x}\right),
\end{equation}
\begin{equation}
F_2\big(a;b,b';c,c';x,y\big)=(1-y)^{-a}
F_2\!\left(a;b,c'-b';c,c';\frac x{1-y},\frac y{y-1}\right),
\end{equation}
\begin{equation}
F_2\big(a;b,b';c,c';x,y\big)=(1-x-y)^{-a}
F_2\!\left(a;c-a,c'-b';c,c';{\mbox{\small $
\frac x{x+y-1},\frac y{x+y-1}$}}\right).
\end{equation}

Also {\em quadratic transformations} are known for Appell functions, see
B.C.~Carlson~\cite{SchlosserCarlson:76}.

\subsection{Reduction formulae}

The transformations of Subsection~\ref{subsec:tf}
readily imply the following reduction formulae
(typically a double series being reduced to a single series):

\smallskip
\noindent$\bullet$\quad $y=x$ in $F_1$:
\begin{subequations}
\begin{equation}
F_1\big(a;b,b';c;x,x\big)=(1-x)^{c-a-b-b'}
\hyp 21{c-a,\,c-b-b'}cx.
\end{equation}
By Euler's transformation this is
\begin{equation}
F_1\big(a;b,b';c;x,x\big)=
\hyp 21{a,\,b+b'}cx.
\end{equation}
\end{subequations}

\noindent$\bullet$\quad $c=b+b'$ in $F_1$:
\begin{equation}\label{eqF12F1}
F_1\big(a;b,b';b+b';x,y\big)=(1-y)^{-a}
\hyp 21{a,\,b}{b+b'}{\frac{x-y}{1-y}}.
\end{equation}

\noindent$\bullet$\quad $c=b$ in $F_2$:
\begin{equation}
F_2\big(a;b,b';b,c';x,y\big)=(1-x)^{-a}
\hyp 21{a,\,b'}{c'}{\frac y{1-x}}.
\end{equation}

\noindent$\bullet$\quad $y=1$ in $F_1$:

Since
$$
F_1\big(a;b,b';c;x,y\big)=
\sum_{m\ge 0}\frac{(a)_m\,(b)_m}{m!\,(c)_m}x^m
\hyp 21{a+m,\,b'}{c+m}y
$$
and
$$
\hyp 21{a,\,b}c1=\Ga\!\begin{bmatrix}c,\,c-a-b\\c-a,\,c-b\end{bmatrix},
\qquad\quad\Re(c-a-b)>0,
$$
we have
\begin{equation}
F_1\big(a;b,b';c;x,1\big)=
\Ga\!\begin{bmatrix}c,\,c-a-b'\\c-a,\,c-b'\end{bmatrix}
\hyp 21{a,\,b}{c-b'}x,
\end{equation}
for $\Re(c-a-b')>0$.

\smallskip
\noindent$\bullet$\quad An $F_1\leftrightarrow F_3$ transformation:

Since
$$
F_1\big(a;b,b';c;x,y\big)=
\sum_{m\ge 0}\frac{(a)_m\,(b)_m}{m!\,(c)_m}x^m
\hyp 21{a+m,\,b'}{c+m}y
$$
and
$$
\hyp 21{a,\,b}cy=(1-y)^{-b}\hyp 21{c-a,\,b}c{\frac y{y-1}},
$$
we have
\begin{align}
F_1\big(a;b,b';c;x,y\big)&=
(1-y)^{-b'}\sum_{m\ge 0}\frac{(a)_m\,(b)_m}{m!\,(c)_m}x^m
\hyp 21{c-a,\,b'}{c+m}{\frac y{y-1}}\notag\\
&=(1-y)^{-b'}F_3\!\left(a,c-a;b,b';c;x,\frac y{y-1}\right).
\end{align}
Hence, any $F_1$ function can be expressed in terms of an $F_3$ function.
The converse is only true when $c=a+a'$.

\smallskip
\noindent$\bullet$\quad $a'=c-a$ and $b'=c-b$ in $F_3$:

Since by Equation \eqref{eqF12F1} the $F_1$ function reduces to an
ordinary $_2F_1$ function when $c=b+b'$, we have
\begin{equation}
F_3\!\left(a,c-a;b,c-b;c;x,\frac y{y-1}\right)=
{\mbox{\small $(1-x)^{-a}(1-y)^{c-b}$}}
\hyp 21{a,\,c-b}c{\frac{y-x}{1-x}}.
\end{equation}

\noindent$\bullet$\quad $c'=a$ in $F_2$:

\begin{equation}
F_2\big(a;b,b';c,a;x,y\big)=(1-y)^{-b'}
F_1\!\left(b;a-b',b';c;x,\frac x{1-y}\right).
\end{equation}
Conversely, any $F_1$ function can be expressed in terms of an
$F_2$ function where $c'=a$.

If further $c=a$, then
\begin{equation}
F_2\big(a;b,b';a,a;x,y\big)=(1-x)^{-b}(1-y)^{-b'}
\hyp 21{b,\,b'}a{\frac{xy}{(1-x)(1-y)}}.
\end{equation}

\subsection{An expansion of an $F_4$ series}

In 1940 and 1941, J.L.~Burchnall and 
T.W.~Chaundy~\cite{SchlosserBC:40,SchlosserBC:41}
gave the following expansion of an $F_4$ series in terms
of products of two hypergeometric $_2F_1$ series:
\begin{multline}
F_4\big(a;b;c,c';x(1-y),y(1-x)\big)\\=
\sum_{m\ge0}\frac{(a)_m\,(b)_m\,(1+a+b-c-c')_m}{m!\,(c)_m\,(c')_m}x^my^m\\
\times\hyp 21{a+m,\,b+m}{c+m}x\hyp 21{a+m,\,b+m}{c'+m}y.
\end{multline}
This expansion has applications to classical orthogonal polynomials.
It can also be used to deduce the double integral representation for $F_4$.
Various special cases are interesting enough to state separately:

\smallskip
\noindent$\bullet$\quad $c'=1+a+b-c$ in $F_4$:

We have the product formula
\begin{equation}
F_4\big(a;b;c,1+a+b-c;x(1-y),y(1-x)\big)=
\hyp 21{a,\,b}cx\hyp 21{a,\,b}{c'}y.
\end{equation}

\noindent$\bullet$\quad $c'=b$ in $F_4$:

Here we have the reduction formula
\begin{align}\label{idc'bF4}
&F_4\big(a;b;c,b;x(1-y),y(1-x)\big)\notag\\&=(1-x)^{-a}(1-y)^{-a}\,
F_1\!\left(a;1+a-c,c-b;c;\frac{xy}{(1-x)(1-y)},\frac x{x-1}\right).
\end{align}

\noindent$\bullet$\quad $c'=b$ and $c=a$ in $F_4$:

Further specialization of \eqref{idc'bF4} gives the quite attractive
summation formula
\begin{subequations}
\begin{equation}
F_4\big(a;b;a,b;x(1-y),y(1-x)\big)=
(1-x)^{1-b}(1-y)^{1-a}(1-x-y)^{-1}.
\end{equation}
Written out in explicit terms, this is
\begin{equation}
\sum_{m\ge0}\sum_{n\ge0}\frac{(a)_{m+n}\,(b)_{m+n}}{m!\,n!\,(a)_m\,(b)_n}
x^m(1-y)^my^n(1-x)^n
=\frac{(1-x)^{1-b}(1-y)^{1-a}}{(1-x-y)}.
\end{equation}
\end{subequations}
For $y=0$ this reduces to I.~Newton's binomial expansion formula
$$
\hyp 10b-x=(1-x)^{-b}.
$$

\section{Related series and extensions of Appell series}

\subsection{Horn functions}\label{subsec:Horn}

In 1931, Jacob Horn~\cite{SchlosserHorn:31} studied
convergent bivariate hypergeometric functions
$\sum_{m,n}f_{m,n}x^my^n$ with certain (degree and other) restrictions
on the two ratios of consecutive terms
$$
\frac{f_{m+1,n}}{f_{m,n}},\qquad\qquad\frac{f_{m,n+1}}{f_{m,n}}.
$$
He arrived at a complete set of $34$ different functions
among which are the Appell functions $F_1,F_2,F_3,F_4$.

They include series such as
\begin{equation}
G_1(a,b,b';x,y):=\sum_{m\ge 0}\sum_{n\ge 0}
\frac{(a)_{m+n}(b)_{n-m}(b')_{m-n}}{m!n!}x^my^n,
\end{equation}
\begin{equation}
H_3(a,b,c;x,y):=\sum_{m\ge 0}\sum_{n\ge 0}
\frac{(a)_{2m+n}(b)_n}{(c)_{m+n}\,m!n!}x^my^n,
\end{equation}
and
\begin{equation}
H_7(a,b,b',c;x,y):=\sum_{m\ge 0}\sum_{n\ge 0}
\frac{(a)_{2m-n}(b)_n(b')_n}{(c)_m\,m!n!}x^my^n.
\end{equation}

\subsection{Kamp\'e de F\'eriet series}
In 1937, J.~Kamp\'e de F\'eriet~\cite{SchlosserKampedeFeriet:37}
introduced the following bivariate extension of the
generalized hypergeometric series:
\begin{multline}
F^{p:q}_{r:s}\left(\begin{matrix}a_1,\dots,a_p:b_1,b_1';\dots;b_q,b_q';\\
c_1,\dots,c_r:d_1,d_1';\dots;d_s,d_s';\end{matrix}\,x,y\right)\\=
\sum_{m\ge0}\sum_{n\ge0}\frac{(a_1)_{m+n}\dots(a_p)_{m+n}}
{(c_1)_{m+n}\dots(c_r)_{m+n}}\frac{(b_1)_m(b_1')_n\dots(b_q)_m(b_q')_n}
{(d_1)_m(d_1')_n\dots(d_s)_m(d_s')_n}\frac{x^my^n}{m!n!}.
\end{multline}

Numerous identities exist for special instances of such series.
For illustration, we list three summation formulae.

\smallskip
\noindent$\bullet$\quad P.W.~Karlsson~\cite{SchlosserKarlsson:94}, 1994: 
\begin{equation}
F^{0:3}_{1:1}\!\left(\begin{matrix}-:a,d-a;b,d-b;c,-c;\\
d:e,d+e-a-b-c;\end{matrix}\,1,1\right)=
\Ga\!\begin{bmatrix}e,\,e+d-a-b-c\\e-c,\,e+d-a-b\end{bmatrix},
\end{equation}
where $\Re(e)>0$ and $\Re(d+e-a-b-c)>0$.

\smallskip
\noindent$\bullet$\quad S.N.~Pitre and
J.~Van der Jeugt~\cite{SchlosserPitreJeugt:96}, 1996: 
\begin{equation}
F^{0:3}_{1:1}\!\left(\begin{matrix}-:a,d-a;b,d-b;c,d-c;\\
d:e,d+e-a-b-c;\end{matrix}\,1,1\right)=
\Ga\!\begin{bmatrix}e,e+d-a-b-c,e-d\\e-a,\,e-b,\,e-c\end{bmatrix}\!,
\end{equation}
where $\Re(e-d)>0$ and $\Re(d+e-a-b-c)>0$. Further
\begin{multline}
F^{0:3}_{1:1}\!\left(\begin{matrix}-:a,d-a;b,d-b;c,e-c-1;\\
d:e,d+e-a-b-c;\end{matrix}\,1,1\right)\\=
\Ga\!\begin{bmatrix}1-a,\,1-b,\,e,\,e-d,\,d+e-a-b-c\\
1-d,\,e-a,\,e-b,\,e-c,\,1+d-a-b\end{bmatrix},
\end{multline}
where $\Re(d+e-a-b-c)>0$, and $d-a$ or $d-b$ is a negative integer.

\subsection{Lauricella series}
In 1893, G.~Lauricella~\cite{SchlosserLauricella:93}
investigated properties of the following four series
$F_A^{(n)}$, $F_B^{(n)}$, $F_C^{(n)}$, $F_D^{(n)}$, of $n$ variables:

\begin{align}
F_A^{(n)}&\big(a;b_1,\dots,b_n;c_1,\dots,c_n;x_1,\dots,x_n\big)
\notag\\
&=\sum_{m_1\ge 0}\dots\sum_{m_n\ge0}\frac{(a)_{m_1+\dots+m_n}\,
(b_1)_{m_1}\dots(b_n)_{m_n}}{(c_1)_{m_1}\dots(c_n)_{m_n}\,
m_1!\dots m_n!}x_1^{m_1}\dots x_n^{m_n},
\end{align}
where $|x_1|+\dots+|x_n|<1$.

\begin{align}
F_B^{(n)}&\big(a_1,\dots,a_n;b_1,\dots,b_n;c;x_1,\dots,x_n\big)
\notag\\
&=\sum_{m_1\ge 0}\dots\sum_{m_n\ge0}\frac{(a_1)_{m_1}\dots(a_n)_{m_n}\,
(b_1)_{m_1}\dots(b_n)_{m_n}}{(c)_{m_1+\dots+m_n}\,
m_1!\dots m_n!}x_1^{m_1}\dots x_n^{m_n},
\end{align}
where $|x_1|,\dots,|x_n|<1$.

\begin{align}
F_C^{(n)}&\big(a;b;c_1,\dots,c_n;x_1,\dots,x_n\big)\notag\\
&=\sum_{m_1\ge 0}\dots\sum_{m_n\ge0}\frac{(a)_{m_1+\dots+m_n}\,
(b)_{m_1+\dots+m_n}}{(c_1)_{m_1}\dots(c_n)_{m_n}\,
m_1!\dots m_n!}x_1^{m_1}\dots x_n^{m_n},
\end{align}
where $|x_1|^{\frac 12}+\dots+|x_n|^{\frac 12}<1$.

\begin{align}
F_D^{(n)}&\big(a;b_1,\dots,b_n;c;x_1,\dots,x_n\big)\notag\\
&=\sum_{m_1\ge 0}\dots\sum_{m_n\ge0}\frac{(a)_{m_1+\dots+m_n}\,
(b_1)_{m_1}\dots(b_n)_{m_n}}{(c)_{m_1+\dots+m_n}\,
m_1!\dots m_n!}x_1^{m_1}\dots x_n^{m_n},
\end{align}
where $|x_1|,\dots,|x_n|<1$.

Certainly, we have 
$$
F_A^{(2)}=F_2,\qquad F_B^{(2)}=F_3,\qquad
F_C^{(2)}=F_4,\qquad F_D^{(2)}=F_1.
$$

Many properties for Lauricella functions, such as
integral representations and partial differential equations,
are given by Appell and Kamp\'e de F\'eriet~\cite{SchlosserAppellKdF:26}.
From the vast amount of material, we single out the following
integral representation of the
Lauricella $F_D^{(n)}$ series as a specific example.

\subsubsection{Integral representation of $F_D^{(n)}$}

The formula
\begin{align}
F_D^{(n)}&\big(a;b_1,\dots,b_n;c;x_1,\dots,x_n\big)\notag\\
&=\Ga\!\begin{bmatrix}c\\a,c-a\end{bmatrix}
\int_0^1 u^{a-1}(1-u)^{c-a-1}(1-ux_1)^{-b_1}\dots(1-ux_n)^{-b_n}\rd u,
\end{align}
where $\Re c>\Re a>0$, is very useful for deriving relations
for $F_D$ series. It
can be easily verified by Taylor expansion of the integrand,
followed by termwise integration.

\subsubsection{Group theoretic interpretations}

A group theoretic interpretation of the Lauricella $F_A^{(n)}$ functions
corresponding to the most degenerate principal series representations
of $\mathrm S\mathrm L(n,\mathbb R)$ was given by
N.Ja.~Vilenkin~\cite{SchlosserVilenkin:70} (see also
\cite[Sec.~16.3.4]{SchlosserVilenkinKlimyk:92}).
Similarly, W.~Miller, Jr.~\cite{SchlosserMiller:72} has shown that
the Lauricella $F_D^{(n)}$ functions transform as
basis vectors corresponding to irreducible representations of
the Lie algebra $sl(n+3,\mathbb C)$ (by which he generalized his previous
observation in \cite{SchlosserMiller:73}
for the $n=2$ case, corresponding to the Appell functions $F_1$).


\bigskip
\textbf{Acknowledgement.}
I would like to thank Tom Koornwinder for pointing out references on
group theoretic interpretations of Lauricella series.

\end{document}